\documentclass[amstex,12pt,russian,amssymb]{article}

\usepackage{mathtext}
\usepackage[cp1251]{inputenc}
\usepackage[T2A]{fontenc}
\usepackage[russian]{babel}
\usepackage[dvips]{graphicx}
\usepackage{amsmath}
\usepackage{amssymb}
\usepackage{amsxtra}
\usepackage{latexsym}
\usepackage{ifthen}

\begin{document}

\newcounter{lemma}
\newcommand{\lemma}{\par \refstepcounter{lemma}%
{\bf Лема \arabic{lemma}.}}

\newcounter{corollary}
\newcommand{\corollary}{\par \refstepcounter{corollary}%
{\bf Наслідок \arabic{corollary}.}}

\newcounter{remark}
\newcommand{\remark}{\par \refstepcounter{remark}%
{\bf Зауваження \arabic{remark}.}}

\newcounter{theorem}
\newcommand{\theorem}{\par \refstepcounter{theorem}%
{\bf Теорема \arabic{theorem}.}}

\newcounter{proposition}
\newcommand{\proposition}{\par \refstepcounter{proposition}%
{\bf Твердження \arabic{proposition}.}}

\newcounter{example}
\newcommand{\example}{\par \refstepcounter{example}%
{\bf Приклад~\arabic{example}.}}

\renewcommand{\refname}{\centerline{\bf Перелік використаних джерел }}

\newcommand{\proof}{{\it Доведення.\,\,}}

\noindent УДК 517.5

{\bf Е.А.~Севостьянов} (Житомирский государственный университет
имени Ивана Франко; Институт прикладной математики и механики НАН
Украины, г.~Славянск)

{\bf С.А.~Скворцов} (Житомирский государственный университет имени
Ивана Франко)

\medskip\medskip
{\bf Є.О.~Севостьянов} (Житомирський державний університет імені
Івана Фран\-ка; Інститут прикладної математики і механіки НАН
України, м.~Слов'янськ)

{\bf С.О.~Скворцов} (Житомирський державний університет імені Івана
Фран\-ка)

\medskip\medskip
{\bf E.A.~Sevost'yanov} (Zhytomyr Ivan Franko State University;
Institute of Applied Ma\-the\-ma\-tics and Mechanics of NAS of
Ukraine, Slov'yans'k)

{\bf S.O.~Skvortsov} (Zhytomyr Ivan Franko State University)

\medskip
{\bf О равностепенной непрерывности семейств отображений,
фиксирующих точку области}

{\bf Про одностайну неперервність сімей відображень, що мають
фіксовану точку області}

{\bf On equicontinuity of families of mappings with a fixed point of
a domain}

\medskip\medskip
Изучено поведение одного класса отображений области евклидова
пространства. Установлено, что указанный класс равностепенно
непрерывен как во внутренних, так и в граничных точках области, если
входящие в него отображения удовлетворяют общему условию нормировки,
а соответствующая характеристика квазиконформности имеет слабый
рост.

\medskip\medskip
Досліджено поведінку одного класу відображень області евклідового
простору. Встановлено, що вказаний клас є одностайно неперервним як
у внутрішніх, так і у межових точках області, якщо відображення, які
входять до нього, задовольняють загальну умову нормування, а
відповідна характеристика квазіконформності має слабке зростання.

\medskip\medskip
The behavior of a class of mappings of a domain of Euclidean space
is studied. It is established that the indicated class is
equicontinuous both at the inner and at the boundary points of the
domain if the mappings contained in it satisfy the general
normalization condition, and the corresponding characteristic of
quasiconformality has a weak growth.

\newpage
{\bf 1. Вступ.} Дану статтю присвячено вивченню відображень з
обмеженим і скінченним спотворенням (див., напр., \cite{MRSY},
\cite{GRY}). Зокрема, в даній замітці ми продовжуємо і дещо
доповнюємо дослідження, проведені в~\cite{SevSkv}.

\medskip
Нещодавно ми дослідили питання щодо локальної і глобальної поведінки
відображень областей, коли їх відповідні образи є змінними
(\cite{SevSkv}). При цьому, один важливий випадок не було враховано,
а саме, не розглянуто ситуацію, коли відображення мають загальну
фіксовану точку. В даній замітці ми покажемо, що відповідні сім'ї
відображень одностайно неперервні всередині і на межі заданої
області за вказаної умови нормування.

\medskip Тут і надалі
\begin{equation}\label{eq49***}
A(x_0, r_1,r_2): =\left\{ x\,\in\,{\Bbb R}^n:
r_1<|x-x_0|<r_2\right\}\,,
\end{equation}
крім того, $M_p(\Gamma)$ позначає $p$-модуль сім'ї кривих $\Gamma$
(див.~\cite{Va}). Нехай $p\geqslant 1$ і $Q:{\Bbb R}^n\rightarrow
[0, \infty]$ -- вимірна за Лебегом функція, що дорівнює нулю зовні
$D.$ Згідно з~\cite[розд.~7.6]{MRSY}, будемо говорити, що
$f:D\rightarrow \overline{{\Bbb R}^n}$ -- {\it кільцеве
$Q$-відображення в точці $x_0\in \overline{D}$ відносно $p$-модуля,}
$x_0\ne \infty,$ якщо при деяких $r_0=r(x_0)>0,$ довільних
$0<r_1<r_2<r_0$ і довільних континуумах $E_1\subset \overline{B(x_0,
r_1)}\cap D,$ $E_2\subset \left(\overline{{\Bbb R}^n}\setminus
B(x_0, r_2)\right)\cap D,$ виконується співвідношення
\begin{equation}\label{eq3*!!}
 M_p\left(f\left(\Gamma\left(E_1,\,E_2,\,D\right)\right)\right)\ \leqslant
\int\limits_{A} Q(x)\cdot \eta^p(|x-x_0|)\ dm(x)\,, \end{equation}
де $\eta : (r_1,r_2)\rightarrow [0,\infty ]$ -- довільна вимірна за
Лебегом функція, що задовольняє нерівність
\begin{equation}\label{eq28*}
\int\limits_{r_1}^{r_2}\eta(r)\, dr \geqslant\ 1\,.
\end{equation}
Відображення $f:D\rightarrow \overline{{\Bbb R}^n}$ називається
кільцевим $Q$-відображенням в $\overline{D}\setminus\{\infty\}$
відносно $p$-модуля, якщо~(\ref{eq3*!!}) виконується для
всіх~$x_0\in \overline{D}\setminus\{\infty\}.$  Дане означення також
можна застосувати до точки~$x_0\in \overline{D}\setminus\{\infty\}$
за допомогою інверсії $\varphi(x)=\frac{x}{|x|^2},$ $\infty\mapsto
0.$ У подальшому $h(x, y)$ позначає хордальну відстань між точками
$x, y\in \overline{{\Bbb R}^n},$ а $h(E)$ -- хордальний діаметр
множини~$E\subset\overline{{\Bbb R}^n}$ (див.~\cite[розд.~1]{MRSY}).

\medskip
Нехай $I$ -- фіксований набір індексів і $D_i,$ $i\in I,$ -- деяка
послідовність областей. Згідно з~\cite[розд.~2.4]{NP}, будемо
говорити, що сім'я областей $\{D_i\}_{i\in I}$ є {\it одностайно
рівномірною відносно $p$-модуля}, якщо для кожного $r>0$ існує число
$\delta>0$ таке, що нерівність
\begin{equation}\label{eq17***}
M_p(\Gamma(F^{\,*},F, D_i))\geqslant \delta
\end{equation}
виконується для всіх~$i\in I$ і довільних континуумів $F, F^*\subset
D_i$ таких, що $h(F)\geqslant r$ і $h(F^{\,*})\geqslant r.$

\medskip
Для $p\geqslant 1,$ заданого числа $\delta>0,$ області $D\subset
{\Bbb R}^n,$ $n\geqslant 2,$ точки $a\in D$ і заданої функції
$Q:D\rightarrow[0, \infty]$ позначимо через $\frak{F}_{Q, a, p,
\delta}(D)$ сім'ю всіх кільцевих $Q$-гомеоморфізмів $f:D\rightarrow
\overline{{\Bbb R}^n}$ в $\overline{D}$ відносно $p$-модуля, які
задовольняють умови $h(f(a), \partial f(D))\geqslant\delta$ і
$h(\overline{{\Bbb R}^n}\setminus f(D))\geqslant \delta.$ Покладемо
$$q_{x_0}(r):=\frac{1}{\omega_{n-1}r^{n-1}}\int\limits_{|x-x_0|=r}Q(x)\,dS\,,$$
де $dS$ -- елемент площі поверхні $S,$ і
$$q^{\,\prime}_{b}(r):=\frac{1}{\omega_{n-1}r^{n-1}}\int\limits_{|x-b|=r}Q^{\,\prime}(x)\,dS\,,$$
$Q^{\,\prime}(x)=\max\{Q(x), 1\}.$ Має місце наступне твердження.

\medskip
\begin{theorem}\label{th3} {\sl\, Припустимо, $p\in (n-1, n],$ область $D$ є локально зв'язною в
кожній точці $x_0\in\partial D$ і області $D_f^{\,\prime}=f(D)$ є
одностайно рівномірними відносно $p$-модуля по всіх~$\frak{F}_{Q, a,
p, \delta}(D).$ Якщо функція $Q$ має скінченне середнє коливання в
$\overline{D},$ або в кожній точці~$x_0\in \overline{D}$ при деякому
$\beta(x_0)>0$ виконується умова
\begin{equation}\label{eq2}
\int\limits_{0}^{\beta(x_0)}\frac{dt}{t^{\frac{n-1}{p-1}}q_{x_0}^{\,\prime\,\frac{1}{p-1}}(t)}=\infty\,,
\end{equation}
то кожне відображення~$f\in\frak{F}_{Q, a, p, \delta}(D)$ має
неперервне продовження в $\overline{D}$ і, крім того, сім'я
$\frak{F}_{Q, a, p, \delta}(\overline{D}),$ що складається з
продовжених таким чином відображень~$\overline{f}:
\overline{D}\rightarrow \overline{{\Bbb R}^n},$ є одностайно
неперервною в~$\overline{D}.$
  }
\end{theorem}

\medskip Як і в~\cite{SevSkv}, ми розглянемо також випадок складних
меж, що відноситься до ситуації простих кінців
(див.~\cite[теореми~3, 4]{SevSkv}). Означення регулярних областей та
областей з локально квазіконформною межею можна знайти, напр.,
в~\cite{KR}, або вказаній публікації~\cite{SevSkv}. Означення
простого кінця, яке тут вживається, та замикання $\overline{D}_P$
області $D$ в термінах простих кінців також можна знайти в
роботі~\cite{KR}. Справедливе наступне твердження.

\medskip
\begin{theorem}\label{th1} {\sl\, Нехай $p\in (n-1, n],$ область $D$ регулярна, а області
$D_f^{\,\prime}=f(D)$ є обмеженими одностайно рівномірними відносно
$p$-модуля по всіх $f\in\frak{F}_{Q, a, p, \delta}(D),$ крім того,
ці області мають локально квазіконформну межу. Якщо функція $Q$ має
скінченне середнє коливання в $\overline{D},$ або в кожній точці
$x_0\in \overline{D}$ при деякому $\beta(x_0)>0$ виконується
умова~(\ref{eq2}), то кожне $f\in\frak{F}_{Q, a, p, \delta}(D)$ має
неперервне продовження~$\overline{f}: \overline{D}_P\rightarrow
\overline{{\Bbb R}^n}$ в $\overline{D}_P$ і, крім того, сім'я
$\frak{F}_{Q, a, p, \delta}(\overline{D})$ усіх продовжених
відображень $\overline{f}: \overline{D}_P\rightarrow \overline{{\Bbb
R}^n}$ є одностайно неперервною в $\overline{D}_P.$
  }
\end{theorem}

\medskip
\begin{remark}\label{rem2}
В теоремі~\ref{th1} одностайну неперервність слід розуміти в сенсі
відображень, що діють між просторами~$(X, d)$ і $\left(X^{\,\prime},
d^{\,\prime}\right),$ де $X=\overline{D}_P$ -- поповнення області
$D$ її простими кінцями, а $d$ -- одна з можливих метрик, що
відповідають топологічному простору~$\overline{D}_P$ (\cite{KR}).
Крім того, $X^{\,\prime}=\overline{{\Bbb R}^n}$ і $d^{\,\prime}$ --
хордальна (сферична) метрика.
\end{remark}

\medskip
{\bf 2. Основні леми.} Наступна лема містить в собі твердження
теореми~\ref{th3} в найбільш загальній ситуації відносно
функції~$Q.$

\medskip
\begin{lemma}\label{lem1}{\sl\, Припустимо, $p\in (n-1, n],$
область $D$ є локально зв'язною в точці $x_0\in\partial D$ і області
$D_f^{\,\prime}=f(D)$ є одностайно рівномірними відносно $p$-модуля
по всіх~$\frak{F}_{Q, a, p, \delta}(D).$ Припустимо, що знайдеться
$\varepsilon_0=\varepsilon_0(x_0)>0$ та вимірна за Лебегом функція
$\psi:(0, \varepsilon_0)\rightarrow [0,\infty]$ така, що
\begin{equation}\label{eq7***} I(\varepsilon,
\varepsilon_0):=\int\limits_{\varepsilon}^{\varepsilon_0}\psi(t)\,dt
< \infty\quad \forall\,\,\varepsilon\in (0, \varepsilon_0)\,,\quad
I(\varepsilon, \varepsilon_0)\rightarrow
\infty\quad\text{при}\quad\varepsilon\rightarrow 0\,,
\end{equation}
і, крім того,
\begin{equation} \label{eq3.7.2}
\int\limits_{A(x_0, \varepsilon, \varepsilon_0)}
Q(x)\cdot\psi^{\,p}(|x-x_0|)\,dm(x) = o(I^p(\varepsilon,
\varepsilon_0))\,,\end{equation}
при $\varepsilon\rightarrow 0$ (де $A(x_0, \varepsilon,
\varepsilon_0)$ визначена в~(\ref{eq49***})).

Тоді кожне відображення~$f\in\frak{F}_{Q, a, p, \delta}(D)$ має
неперервне продовження в точку $x_0$ і, крім того, сім'я
$\frak{F}_{Q, a, p, \delta}(D\cup\{x_0\}),$ що складається з
продовжених таким чином відображень~$\overline{f}:
D\cup\{x_0\}\rightarrow \overline{{\Bbb R}^n},$ є одностайно
неперервною в~$D\cup\{x_0\}.$
 }
\end{lemma}

\medskip
\begin{proof}
Одностайна неперервність сім'ї $\frak{F}_{Q, a, p,
\delta}(D\cup\{x_0\})$ всередині області $D,$ а також наявність
неперервного продовження кожного $f\in\frak{F}_{Q, a, p,
\delta}(D\cup\{x_0\})$ випливають аналогічно до міркувань, наведених
на початку доведення~\cite[лема~1]{SevSkv}.

\medskip
Залишилось довести одностайну неперервність $\frak{F}_{Q, a, p,
\delta}(D\cup\{x_0\})$ в точці $x_0.$ Доведемо це від супротивного.
Припустимо, що сім'я відображень~$\frak{F}_{Q, a, p,
\delta}(D\cup\{x_0\})$ не є одностайно неперервною в точці $x_0.$
Тоді існує $\varepsilon_*>0$ з наступною умовою: знайдеться
послідовність $x_m\in D,$ $x_m\rightarrow x_0$ при
$m\rightarrow\infty,$ і $f_m\in \frak{F}_{Q, a, p, \delta}(D)$ такі,
що
\begin{equation}\label{eq5}
h(f_m(x_0), f_m(x_m))\geqslant \varepsilon_*.
\end{equation}
Оскільки $f_m$ мають неперервне продовження в точку $x_0,$ з
урахуванням~(\ref{eq5}) випливає існування послідовності
$x^{\,\prime}_m\in D,$ $x^{\,\prime}_m\rightarrow x_0$ при
$m\rightarrow\infty,$ такої, що
\begin{equation}\label{eq6}
h(f_m(x^{\,\prime}_m), f_m(x_m))\geqslant \varepsilon_*/2.
\end{equation}
Оскільки $D$ є локально зв'язною на своїй межі, то ми можемо
вважати, що точки $x_m$ і $x^{\prime}_m$ належать околу $V_m$ точки
$x_0$ такому, що $W_m:=V_m\cap D$ зв'язна і $W_m\subset B(x_0,
2^{\,-m}).$ З'єднаємо точки $x_m$ і $x^{\prime}_m$ кривою
$\gamma_m\subset W_m\subset B(x_0, 2^{\,-m}).$

\medskip
Застосуємо тепер наступну конструкцію. Якщо $\partial D$ містить
принаймні одну скінченну точку, крім $x_0,$ позначимо її символом
$y_0.$ В протилежному випадку, покладемо $y_0:=\infty.$ Оскільки
область $D$ є локально зв'язною в усіх точках межі, точку $a$ можна
з'єднати кривою з точкою $y_0,$ яка цілком належить $D,$ крім своєї
кінцевої точки $y_0;$ див. з цього
приводу~\cite[пропозиція~13.2]{MRSY}. Останню криву позначимо через
$E.$ Без обмеження загальності, можна вважати, що $E\subset
D\setminus B(x_0, \varepsilon_0),$ де $\varepsilon_0$ -- число з
умови леми. Зауважимо, що $C(y_0, f_m)\subset
\partial f_m(D),$ див.,
напр.,~\cite[пропозиція~13.5]{MRSY}. Доведемо тепер існування точки
$a_m\in E,$ такої що
\begin{equation}\label{eq46*}
h(f_m(a), f_m(a_m))\geqslant (1/2)\cdot h (f_m(a),
\partial f_m(D))\geqslant\delta/2\,.
\end{equation}
Дійсно, зафіксуємо $m=1,2,\ldots .$ Нехай $z_k\in E,$ $k=1,2\ldots $
-- довільна послідовність, для якої $z_k\rightarrow y_0$ при
$k\rightarrow\infty.$ Оскільки простір $\overline{{\Bbb R}^n}$ є
компактним, то можна вважати, що послідовність $f_m(z_k)$ також
збігається до деякої точки $y_0^{\,\prime}$ по хордальній метриці
при $k\rightarrow\infty.$ Оскільки при гомеоморфізмах $C(y_0,
f_m)\subset
\partial f_m(D)$ (див., напр.,~\cite[пропозиція~13.5]{MRSY}),
то $y_0^{\,\prime}\in\partial f_m(D).$ Оскільки $f_m(z_k)\rightarrow
y_0^{\,\prime}$ при $k\rightarrow\infty,$ то для числа $\delta/2$
знайдеться $k_0\in {\Bbb N}$ такий, що
\begin{equation}\label{eq1D}
h(f_m(z_k), y_0^{\,\prime})<\delta/2
\end{equation}
при всіх $k\geqslant k_0.$
Покладемо $a_m:=z_{k_0}.$ Тоді, за нерівністю трикутника,
співвідношення~(\ref{eq1D}), а також за означенням класу
відображень~$\frak{F}_{Q, a, p, \delta}(D),$ маємо:
$$\delta\leqslant h(f_m(a), \partial f_m(D))\leqslant$$
\begin{equation}\label{eq1E}
\leqslant h(f_m(a), y^{\,\prime}_0)\leqslant h(f_m(a), f_m(a_m))+
h(f_m(a_m), y^{\,\prime}_0) \leqslant
\end{equation}
$$\leqslant  h(f_m(a), f_m(a_m))+ \delta/2\,,$$
або, переносячи $\delta/2$ в ліву частину~(\ref{eq1E}),
$$\delta/2\leqslant h(f_m(a), \partial f_m(D))\leqslant  h(f_m(a), f_m(a_m))\,.$$
Останнє співвідношення доводить~(\ref{eq46*}).

\medskip
Замкнену підкриву кривої $E,$ що з'єднує точки $a_0$ і $a$ в $D,$
позначимо через $F_m$ (див. малюнок~\ref{fig1}).
\begin{figure}[h]
\centerline{\includegraphics[scale=0.7]{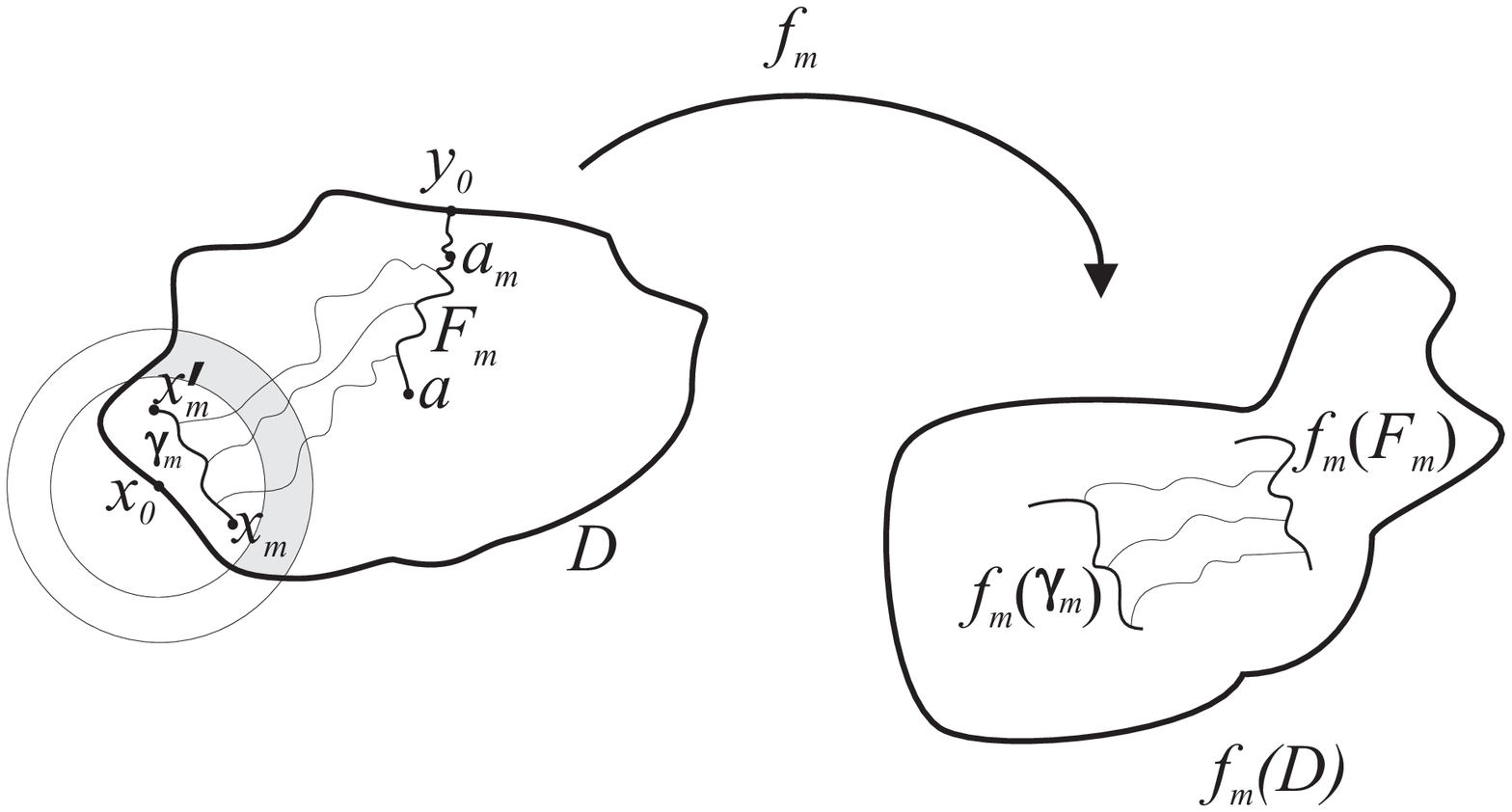}} \caption{До
доведення леми~\ref{lem1}}\label{fig1}
\end{figure}
Нехай $r_0\in (0, \varepsilon_0)$ є таким, що $I(r,
\varepsilon_0)>0$ при $0<r<r_0$ (це можливо з огляду на
умову~$I(\varepsilon, \varepsilon_0)\rightarrow \infty$ при
$\varepsilon\rightarrow 0$). Нехай також $m_0\in {\Bbb N}$ є таким,
що $2^{\,-m}<r_0$ при $m\geqslant m_0.$ Тоді для цих самих $m$
розглянемо сім'ю вимірних функцій
$$\eta_m(t)= \left\{
\begin{array}{rr}
\psi(t)/I(2^{\,-m}, \varepsilon_0), &   t\in (2^{\,-m},\varepsilon_0),\\
0,  &  t\not\in (2^{\,-m}, \varepsilon_0)\,.
\end{array}
\right.$$
Зауважимо, що функції~$\eta_m$ задовольняють
співвідношення~(\ref{eq28*}) при $r_1=2^{\,-m}$ і
$r_2=\varepsilon_0,$ відповідно. В такому випадку, за
співвідношенням~(\ref{eq3*!!}) при $m\rightarrow\infty$
$$M_p(f_m(\Gamma(|\gamma_m|, F_m, D)))=$$
\begin{equation}\label{eq37***}=M_p(\Gamma(f_m(|\gamma_m|), f_m(F_m),
f_m(D)))\leqslant \alpha(2^{\,-m})\rightarrow 0\,,
\end{equation}
де $\alpha(r)$ -- деяка функція, що прямує до нуля при $r\rightarrow
0$ і існування якої обумовлено співвідношенням~(\ref{eq3.7.2}).

\medskip
Проте, оскільки сім'я областей $f_m(D)$ є одностайно рівномірною
відносно $p$-модуля, то враховуючи співвідношення~(\ref{eq46*}),
маємо:
\begin{equation}\label{eq1}
M_p(\Gamma(f_m(|\gamma_m|), f_m(F_m), f_m(D)))\geqslant\delta_*>0
\end{equation}
за певного~$\delta_*>0.$ Співвідношення~(\ref{eq37***}) і
(\ref{eq1}) суперечать одне одному, що і доводить лему.
\end{proof}$\Box$

\medskip
Твердження теореми~\ref{th3} безпосередньо випливає з
леми~\ref{lem1} і \cite[пропозиція~2, деталі доведення
теореми~2]{Sev$_5$}.~$\Box$

\medskip
Сформулюємо і доведемо також лему про аналогічні класи відображень
для випадку, коли область $D$ має межу більш складної структури.

\medskip
\begin{lemma}\label{lem2} {\sl\, Нехай $p\in (n-1, n],$ область $D$ регулярна, а області
$D_f^{\,\prime}=f(D)$ є обмеженими одностайно рівномірними відносно
$p$-модуля по всіх $f\in\frak{F}_{Q, a, p, \delta}(D),$ крім того,
ці області мають локально квазіконформну межу. Припустимо, що
знайдеться $\varepsilon_0=\varepsilon_0(x_0)>0$ та вимірна за
Лебегом функція $\psi:(0, \varepsilon_0)\rightarrow [0,\infty]$
така, що виконуються умови~(\ref{eq7***})--(\ref{eq3.7.2}).

Тоді кожне $f\in\frak{F}_{Q, a, p, \delta}(D)$ має неперервне
продовження~$\overline{f}: \overline{D}_P\rightarrow \overline{{\Bbb
R}^n}$ в $\overline{D}_P$ і, крім того, сім'я $\frak{F}_{Q, a, p,
\delta}(\overline{D})$ усіх продовжених відображень $\overline{f}:
\overline{D}_P\rightarrow \overline{{\Bbb R}^n}$ є одностайно
неперервною в $\overline{D}_P.$
  }
\end{lemma}

\medskip
\begin{proof}
Одностайна неперервність сім'ї $\frak{F}_{Q, a, p,
\delta}(\overline{D})$ всередині області $D,$ а також наявність
неперервного продовження кожного $f\in\frak{F}_{Q, a, p,
\delta}(\overline{D})$ випливають аналогічно до міркувань, наведених
на початку доведення~\cite[лема~1]{SevSkv}.

\medskip
Доведемо одностайну неперервність $\frak{F}_{Q, a, p,
\delta}(\overline{D})$ в $E_D:=\overline{D}_P\setminus D.$
Припустимо протилежне, а саме, що існують $\varepsilon_*>0,$ $P_0\in
E_D,$ послідовність $x_m\in \overline{D}_P,$ $x_m\rightarrow P_0$
при $m\rightarrow\infty,$ і відображення $f_m\in \frak{F}_{Q, a, p,
\delta}(\overline{D})$ такі, що
\begin{equation}\label{eq1A}
h(f_m(x_m), f_m(P_0))\geqslant\varepsilon_*\,,\quad m=1,2,\ldots\,.
\end{equation}
Враховуючи~(\ref{eq1A}) і міркуючи так, як і при доведенні
попередньої леми~\ref{lem1}, ми можемо вважати, що $x_m\in D,$ крім
того, існує ще одна послідовність $x^{\,\prime}_m\in
\overline{D}_P,$ $x^{\,\prime}_m\rightarrow P_0$ при
$m\rightarrow\infty,$ така, що
\begin{equation}\label{eq1*}
h(f_m(x_m), f_m(x^{\,\prime}_m))\geqslant\varepsilon_*/2\,,\quad
m=1,2,\ldots\,.
\end{equation}
Нехай $d_m$ -- послідовність областей розрізів~$\sigma_m$, що
відповідає кінцю $P_0$ така, що $\sigma_m\subset S(x_0, r_m),$
$x_0\in\partial D$ і $r_m\rightarrow 0$ при $m\rightarrow\infty$
(див. \cite[лема~2]{KR}). Без обмеження загальності, можна вважати,
що $x_m, x^{\,\prime}_m\in d_m.$ Нехай тоді $\gamma_m$ -- крива, що
з'єднує точки $x_m$ і $x^{\,\prime}_m$ всередині $d_m.$

\medskip
Оскільки область $D$ регулярна, простір $\overline{D}_P$ містить не
менше двох простих кінців $P_1$ і $P_2\in E_D.$ Нехай $P_1\subset
E_D$ -- простий кінець, що не співпадає з $P_0.$ Припустимо, $G_m,$
$m=1,2,\ldots,$ -- послідовність областей, яка відповідає простому
кінцю $P_1.$ Оскільки при кожному $m=1,2,\ldots $ відображення $f_m$
має неперервне продовження на $\overline{D}_P,$ можна підібрати
послідовність $\zeta_m\in G_m,$ $\zeta_m\rightarrow P_1$ при
$m\rightarrow\infty,$ такою, що $h(f_m(\zeta_m),
f_m(P_1))\rightarrow 0$ при $m\rightarrow\infty.$ Зауважимо, що при
всіх $m\geqslant m_0$ і деякому $m_0\in {\Bbb N}$
\begin{equation}\label{eq3D}
h(f_m(a), f_m(\zeta_m))\geqslant h(f_m(a), f_m(P_1))-h(f_m(\zeta_m),
f_m(P_1))\geqslant \delta/2\,,
\end{equation}
де, як звично, $h(x, y)$ позначає хордальну відстань між $x$ і $y.$
Побудуємо послідовність континуумів $K_m,$ $m=1,2,\ldots,$ наступним
чином. З'єднаємо точки $\zeta_1$ і $a$ довільною кривою в $D,$ котру
позначимо $K_1.$ Далі, з'єднаємо точки $\zeta_2$ і $\zeta_1$ кривою
$K_1^{\prime},$ в $G_1.$ Поєднавши криві $K_1$ і $K_1^{\prime},$
отримаємо криву $K_2,$ що з'єднує точки $a$ і $\zeta_2.$ І так далі.
Нехай на деякому кроці маємо криву $K_m,$ що з'єднує точки $\zeta_m$
і $a.$ З'єднаємо точки $\zeta_{m+1}$ і $\zeta_m$ кривою
$K_m^{\,\prime},$ яка лежить в $G_m.$ Поєднавши між собою криві
$K_m$ і $K_m^{\,\prime},$ отримаємо криву $K_{m+1}.$ І так далі
(див. малюнок~\ref{fig2}).
\begin{figure}[h]
\centerline{\includegraphics[scale=0.7]{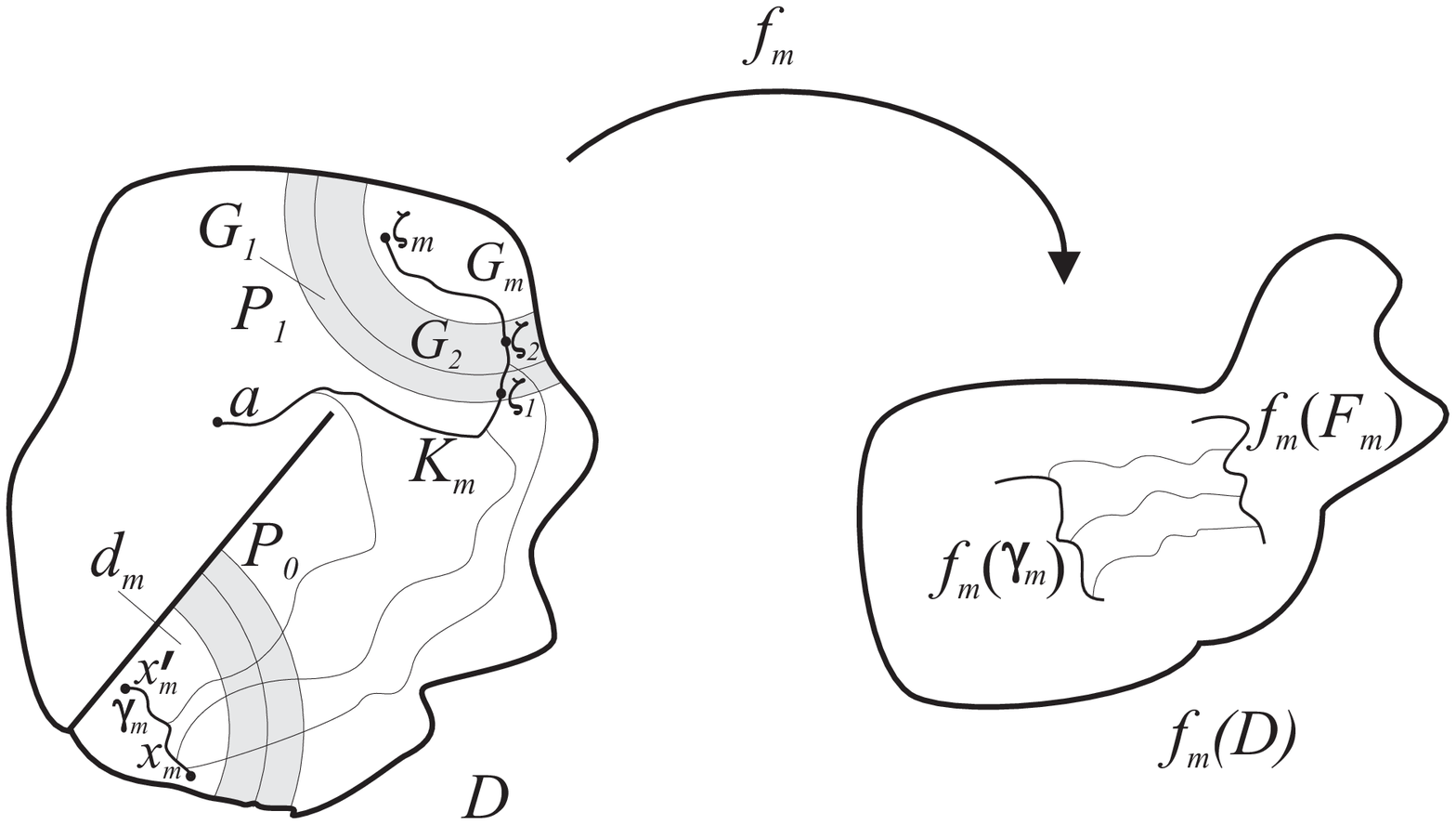}} \caption{До
доведення леми~\ref{lem2}}\label{fig2}
\end{figure}

Покажемо, що знайдеться номер $m_1\in {\Bbb N},$ такий що
\begin{equation}\label{eq4D}
d_m\cap K_m=\varnothing\quad\forall\quad m\geqslant m_1\,.
\end{equation}
Припустимо, що~(\ref{eq4D}) не виконується, тоді знайдеться
зростаюча послідовність номерів~$m_k\rightarrow\infty,$
$k\rightarrow\infty,$ і точок~$\xi_k\in K_{m_k}\cap d_{m_k},$
$m=1,2,\ldots,\,.$ Тоді $\xi_k \rightarrow P_0$ при
$k\rightarrow\infty.$

\medskip
Зауважимо, що можливі два випадки: або всі елементи $\xi_k$ при
$k=1,2,\ldots$ належать $D\setminus G_1,$ або знайдеться номер $k_1$
такий, що $\xi_{k_1}\in G_1.$ Далі, розглянемо
послідовність~$\xi_k,$ $k>k_1.$ Зауважимо, що можливі два випадки:
або $\xi_k$ при $k>k_1$ належать $D\setminus G_2,$ або знайдеться
$k_2>k_1$ такий, що $\xi_{k_2}\in G_2.$ І так далі. Припустимо, що
елемент $\xi_{k_{l-1}}\in G_{l-1}$ вже побудований. Зауважимо, що
можливі два випадки: або $\xi_k$ належать $D\setminus G_l$ при
$k>k_{l-1},$ або знайдеться номер $k_l>k_{l-1}$ такий, що
$\xi_{k_l}\in G_l.$ І т.д. Ця процедура може бути як скінченною, так
і нескінченною, в зв'язку з чим маємо дві можливі ситуації:

1) або знайдуться номери $n_0\in {\Bbb N}$ і $l_0\in {\Bbb N}$ такі,
що $\xi_k\in D\setminus G_{n_0}$ при всіх $k>l_0;$

2) або для кожного $l\in {\Bbb N}$ знайдеться елемент $\xi_{k_l}$
такий, що $\xi_{k_l}\in G_l,$ причому послідовність $k_l$ є
зростаючою по $l\in {\Bbb N}.$

\medskip
Розглянемо кожен з цих випадків окремо і покажемо, що в обох з них
ми приходимо до суперечності. Нехай має місце ситуація~1), тоді
зауважимо, що всі елементи послідовності~$\xi_k$ належать~$K_{n_0},$
звідки випливає існування підпослідовності~$\xi_{k_r},$
$r=1,2,\ldots,$ збіжної при $r\rightarrow\infty$ до деякої точки
$\xi_0\in D.$ Проте, $\xi_k\in d_{m_k}$ і, отже, $\xi_0\in
\bigcap\limits_{m=1}^{\infty} \overline{d_m}\subset
\partial D$ (див.~\cite[пропозиція~1]{KR}).
Отримана суперечність вказує на неможливість випадку~1). Нехай має
місце випадок~2), тоді одночасно $\xi_k\rightarrow P_0$ і
$\xi_k\rightarrow P_1$ при $k\rightarrow\infty.$ Оскільки
простір~$\overline{D}_P$ є метричним (див.~\cite[зауваження~3]{KR},
то звідси, за нерівністю трикутника, випливає, що $P_1=P_0,$ що
суперечить обранню~$P_1.$ Отримана суперечність вказує на
справедливість співвідношення~(\ref{eq4D}).

\medskip Покладемо тепер
$\widetilde{\varepsilon_0}=\min\{\varepsilon_0, r_{m_1+1}\},$ і
нехай $M_0$ -- натуральне число, таке що
$r_m<\widetilde{\varepsilon_0}$ при всіх $m\geqslant M_0.$
Розглянемо сім'ю вимірних функцій
$$\eta_{m}(t)= \left\{
\begin{array}{rr}
\psi(t)/I(r_m, \widetilde{\varepsilon_0}), &   t\in (r_m, \widetilde{\varepsilon_0}),\\
0,  &  t\not\in (r_m, \widetilde{\varepsilon_0})\,,
\end{array}
\right.$$
%
%
де, як і раніше, $I(a, b)$ визначена співвідношенням~$I(a,
b)=\int\limits_a^b\psi(t)\,dt.$ Зауважимо, що функція $\eta_m$ при
кожному фіксованому $m\geqslant M_0$ задовольняє
співвідношення~(\ref{eq28*}) ($r_1:=r_m$ і
$r_2:=\widetilde{\varepsilon_0}$). За співвідношенням~(\ref{eq4D})
та за означенням розрізів~$\sigma_m\subset S(x_0, r_m),$ маємо:
$$\Gamma\left(|\gamma_m|, K_m, D\right)>\Gamma(S(x_0, r_m), S(x_0,
\widetilde{\varepsilon_0}), D)\,.$$
Отже,
$$f_m(\Gamma\left(|\gamma_m|, K_m, D\right))>f_m(\Gamma(S(x_0, r_m), S(x_0,
\widetilde{\varepsilon_0}), D))\,,$$ звідки, за означенням
класу~$\frak{F}_{Q, a, p, \delta}(\overline{D}),$ з урахуванням
співвідношень~(\ref{eq7***})--(\ref{eq3.7.2}), будемо мати, що
\begin{equation}\label{eq39}
M_p(f_m(\Gamma(|\gamma_m|, K_m, D)))\leqslant M_p(f_m(\Gamma(S(x_0,
r_m), S(x_0, \widetilde{\varepsilon_0}), D))\leqslant
\alpha(r_m)\,,\quad m\geqslant M_0\,,
\end{equation}
де $\alpha(r_m)\rightarrow 0,$ $m\rightarrow\infty,$ -- деяка
функція, існування якої зумовлено співвідношенням~(\ref{eq3.7.2}).
Співвідношення~(\ref{eq39}) суперечить одностайній рівномірності
послідовності областей~$D^{\,\prime}_m:=f_m(D).$ Дійсно,
$h(f_m(K_m))\geqslant \delta/2$ згідно з~(\ref{eq3D}), а
$h(f_m(|\gamma_m|))>\varepsilon_*/2$ за
співвідношенням~(\ref{eq1*}). Отже, оскільки послідовність областей
$D^{\,\prime}_m:=f_m(D)$ є рівномірною, для деякого $\delta_*>0$ і
всіх $m=1,2,\ldots, $ маємо:
$$M_p(f_m(\Gamma(|\gamma_m|, K_m, D)))=M_p(\Gamma(f_m(|\gamma_m|), f_m(K_m), f_m(D)))\geqslant
\delta_*>0\,,$$
що суперечить співвідношенню~(\ref{eq39}). Отримана суперечність
вказує на невірність припущення в~(\ref{eq1A}). Лема
доведена.~$\Box$
\end{proof}

\medskip
Твердження теореми~\ref{th1} безпосередньо випливає з
леми~\ref{lem2} і \cite[пропозиція~2, деталі доведення
теореми~2]{Sev$_5$} (див. також~\cite[лема~2.3.1]{Sev$_7$}).~$\Box$

\medskip
\begin{remark}\label{rem1}
В жодному з результатів роботи (теореми~\ref{th3}, \ref{th1};
леми~\ref{lem1}, \ref{lem2}) не можна позбутися умови, що відповідна
сім'я відображень $\frak{F}_{Q, a, p, \delta}(D)$ складається з
гомеоморфізмів. Це, зокрема, підтверджується прикладом аналітичних
функцій $f_n(z)=z^n,$ $n=1,2,\ldots ,$ що відображають одиничний
круг на себе ($Q(z)\equiv 1,$ $D={\Bbb D}=\{z\in {\Bbb C}:
|z|<1\}$). Відображення $f_n$ фіксують точку $z_0=0,$ в той самий
час, сім'я $\{f_n\}_{n=1}^{\infty}$ не є одностайно неперервною на
одиничному колі, що можна перевірити шляхом прямих обчислень.

У випадку, коли відображена область є фіксованою областю
квазіекстремальної довжини, теорема~\ref{th3} встановлена
в~\cite[твердження~1$^{\,\prime}$]{Sev$_2$}, а теорему~\ref{th1} --
для деякого підкласу відображень, що розглядаються в даній статті,
в~\cite[теорема~4]{Sev$_5$}.
\end{remark}

\medskip
\begin{example}\label{ex1}
Розглянемо наступну сім'ю відображень:
$$f_n(z)= \left\{
\begin{array}{rr}z/(2-2^{-(n-1)}), &   |z|<1-2^{-n},\\
z(1-a_n(1-|z|)), & |z|\geqslant 1-2^{-n}\,,
\end{array}
\right.\quad n=1,2,\ldots \,.$$
де $a_n=(2^{n-1}-1)/(1-2^{-n}).$ Зауважимо, що при кожному
фіксованому $n=1,2\ldots$ відображення $f_n^1(z):=z/[2-2^{-(n-1)}]$
є конформним, тому його внутрішня дилатація $K_I(z, f_n^1)$
дорівнює~1. (Означення внутрішньої дилатації $K_I(x, f)$
відображення $f$ у точці $x,$ а також радіальної і тангенсальної
дилатацій $\delta_r(x)$ і $\delta_{\tau}(x)$ можна знайти,
наприклад, у~\cite[співвідношення~(2.30), деталі доведення
пропозиції~6.3]{MRSY}). Крім того, використовуючи підхід і
позначення, застосовані при розгляданні пропозиції~6.3
в~\cite{MRSY}, для відображень $f_n^2(z)=z[1-a_n(1-|z|)]$ при
$1>|z|\geqslant 1-2^{-n}$ маємо:
$$\delta_{\tau}(z)=\frac{|f_n^2(z)|}{|z|}=1-a_n+a_n|z|,$$
$$\delta_r(z)=\left|\frac{\partial |f_n^2(z)|}{\partial |z|}\right|=
1-a_n+2a_n|z|\,.$$
Звідси $\delta_r(z)\geqslant \delta_{\tau}(z),$ отже,
$$K_I(f_n^2, z)=\frac{\delta^{n-1}_{\tau}\delta_r}{\delta^n_{\tau}}=
\frac{(1-a_n+a_n|z|)^{n-1}(1-a_n+2a_n|z|)}{(1-a_n+a_n|z|)^n}=$$
\begin{equation}\label{eq2A}
=\frac{1-a_n+2a_n|z|}{1-a_n+a_n|z|}=1+\frac{a_n|z|}{1-a_n+a_n|z|}\leqslant
\end{equation}
$$\leqslant 1+\frac{a_n}{1-a_n+a_n(1-2^{-n})}=1+(2^{n-1}-1)/2\,.$$
Слід зауважити, що оцінка~(\ref{eq2A}) є точною і досягається в
точках кола $S(0, 1-2^{-n}).$ З~(\ref{eq2A}) випливає, що при
кожному фіксованому~$n\in {\Bbb N}$ внутрішня дилатація
відображення~$f_n^2$ є обмеженою деякою сталою
$c_n:=1+(2^{n-1}-1)/2,$ тому $f_n^2$ є кільцевим $Q_n$-відображенням
при $Q_=c_n$ (див., напр., \cite[теореми~8.1, 8.6]{MRSY}). Більше
того, $f_n^2$ є квазіконформним у кільці~$\{1>|z|> 1-2^{-n}\}$
(див., напр.,~\cite[теорема~34.6]{Va}).

\medskip
За теоремою про усувність аналітичних дуг (див.,
напр.,~\cite[теорема~I.8.3]{LV}) відображення $f_n$ також є
$c_n$-квазіконформним. Отже, $f_n$ є кільцевим $Q_n$-відображенням
при $Q_n=c_n$ (див., напр., \cite[теореми~8.1, 8.6]{MRSY}).

\medskip
Зауважимо, що $f_n(0)=0,$ $n=1,2,\ldots ,$ причому $f_n$
відображають одиничний круг ${\Bbb D}=\{z\in {\Bbb C}: |z|<1\}$ на
себе при кожному $n\in {\Bbb N}.$ В той самий час, $f_n$ збігаються
локально рівномірно до гомеоморфізму $f(z)=z/2,$ що переводить
${\Bbb D}$ на круг $B(0, 1/2)=\{z:|z|<1/2\}$ вдвічі меншого
радіуса~(!). Ми підкреслюємо тут, що локально рівномірна збіжність
автоморфізмів $f_n$ одиничного круга не гарантує, що граничне
відображення також буде автоморфізмом ${\Bbb D}$ (важливо, що це не
так, навіть коли граничне відображення $f$ є гомеоморфізмом). Можна
довести, що побудована послідовність відображень $f_n$ не є
рівномірно збіжною до $f$ в одиничному крузі. Тобто, їхня збіжність
є локально рівномірною, але не рівномірною.

\medskip
Зауважимо, що елементи послідовності~$f_n$ не мають спільної
мажоранти $Q$ у нерівності~(\ref{eq3*!!}), $D:={\Bbb D},$ такої, яка
б задовольняла принаймні одну з можливих умов $FMO,$ (\ref{eq2}),
(\ref{eq7***}) чи~(\ref{eq3.7.2}). Для того, щоб переконатися в
останньому, припустимо протилежне. Тоді, очевидно,
$f_n\in\frak{F}_{Q, 0, n, \delta}(D)$ при кожному $n=1,2,\ldots ,$
де $\delta=\min\{h(0,
\partial {\Bbb D}), h(\overline{C}\setminus{\Bbb D})\}.$ Зауважимо
також, що область ${\Bbb D}$ є рівномірною згідно теореми  Няккі, як
плоска область, що є локально зв'язною на своїй межі і має одну
межову компоненту (див.~\cite[теорема~6.2 і наслідок~6.8]{Na$_3$}).
Тоді за теоремою~\ref{th3}, або лемою~\ref{lem1}, відповідно, сім'я
відображень $f_n,$ $n=1,2,\ldots ,$ є одностайно неперервною
в~$\partial {\Bbb D}.$ Оскільки хордальна і евклідова метрики є
еквівалентними на компактах в~$\overline{{\Bbb R}^n},$ ми можемо
розуміти цю одностайну неперервність у звичайному, евклідовому
сенсі.

\medskip
З іншого боку, зафіксуємо, наприклад, точку $z_0:=1\in
\partial {\Bbb D}.$ Маємо: $f_n(1)=1.$ Розглянемо послідовність
$z_k=1-1/k,$ $k=1,2,\ldots \,.$ Оскільки послідовність $f_n(z)$
збігається локально рівномірно до відображення $f(z)=z/2,$ для
кожного фіксованого $k=1,2,\ldots$ знайдеться номер $n_k\in {\Bbb
N}$ такий, що
$$|f_{n_k}(z_k)-z_k/2|<1/k\,,\quad k=1,2,\ldots \,.$$
Можна вважати, що послідовність номерів $n_k$ є зростаючою, тобто,
$n_1<n_2<n_3<\ldots  .$ Тоді з останньої нерівності, за нерівністю
трикутника, маємо, що
$$|f_{n_k}(z_k)-1|\geqslant |z_k/2-1|-|f_{n_k}(z_k)-z_k/2|\geqslant 1/2-1/k\geqslant 1/4$$
для всіх $k\geqslant 4.$ Остання нерівність суперечить одностайній
неперервності сім'ї відображень $f_n$ в точці~1. Отримана
суперечність спростовує припущення про те, що $f_n$ мають спільну
мажоранту $Q$ у нерівності~(\ref{eq3*!!}), $D:={\Bbb D},$ яка
задовольняє принаймні одну з умов $FMO,$ (\ref{eq2}),
(\ref{eq7***}), або~(\ref{eq3.7.2}).
\end{example}

КОНТАКТНАЯ ИНФОРМАЦИЯ

\medskip
\noindent{{\bf Евгений Александрович Севостьянов} \\
{\bf 1.} Житомирский государственный университет им.\ И.~Франко\\
кафедра математического анализа, ул. Большая Бердичевская, 40 \\
г.~Житомир, Украина, 10 008 \\
{\bf 2.} Институт прикладной математики и механики
НАН Украины, \\
отдел теории функций, ул.~Добровольского, 1 \\
г.~Славянск, Украина, 84 100\\
e-mail: esevostyanov2009@gmail.com}

\medskip
\noindent{{\bf Сергей Александрович Скворцов} \\
Житомирский государственный университет им.\ И.~Франко\\
кафедра математического анализа, ул. Большая Бердичевская, 40 \\
г.~Житомир, Украина, 10 008 \\ e-mail: serezha.skv@gmail.com }

\end{document}